\documentstyle{article}

\newcommand{\beqa}{\begin{eqnarray*}}
\newcommand{\eeqa}{\end{eqnarray*}}
\newcommand{\beqn}{\begin{eqnarray}}
\newcommand{\eeqn}{\end{eqnarray}}
\newcommand{\bm}[1]{\mbox{\boldmath ${#1}$}}
\newcommand{\R}{\mbox{$I\!\!R$}}

\newcommand{\alt}[4]{\left\{ \begin{array}{ll}#1 & \mbox{ if
	\quad}#2 \\ #3 & \mbox{ if \quad}#4 \end{array} \right.}
\newcommand{\alter}[6]{\left\{ \begin{array}{ll}#1 & \mbox{ if
	\quad}#2 \\ #3 & \mbox{ if \quad} #4 \\ #5 & \mbox{ if \quad}#6
		\end{array} \right.}

\newcommand{\ra}{\rightarrow}
\newcommand{\lgra}{\longrightarrow}
\newcommand{\Lglra}{\Longleftrightarrow}

\renewcommand{\a}{\alpha}
\newcommand{\e}{\varepsilon}
\newcommand{\p}{\pi}
\newcommand{\ph}{\phi}

\newcommand{\rec}[1]{\frac{1}{#1}}
\newcommand{\f}{\frac}
\newcommand{\iy}{\infty}
\newcommand{\qed}{\hfill\rule{3mm}{3mm}}
\newcommand{\ds}{\displaystyle}

\newcounter{cnt1}
\newcounter{cnt2}
\newcounter{cnt3}
\newcommand{\blr}{\begin{list}{$($\roman{cnt1}$)$}
	{\usecounter{cnt1} \setlength{\topsep}{0pt}
		\setlength{\itemsep}{0pt}}}
\newcommand{\bla}{\begin{list}{$($\alph{cnt2}$)$}
	{\usecounter{cnt2} \setlength{\topsep}{0pt}
		\setlength{\itemsep}{0pt}}}
\newcommand{\bln}{\begin{list}{$($\arabic{cnt3}$)$}
	{\usecounter{cnt3} \setlength{\topsep}{0pt}
		\setlength{\itemsep}{0pt}}}
\newcommand{\el}{\end{list}}

\newtheorem{thm}{Theorem}[section]
\newtheorem{lem}[thm]{Lemma}
\newtheorem{cor}[thm]{Corollary}

\title{EXTREME CONTRACTIONS IN ${\cal L}(\ell^p_2, \ell^q_2)$
AND THE MAZUR INTERSECTION PROPERTY IN $\ell^p_2 \otimes_{\p}
\ell^q_2$}
\author{Pradipta Bandyopadhyay \and A.\ K.\ Roy\\ Stat--Math 
Division, Indian Statistical Institute\\ 203, B.\ T.\ Road, 
Calcutta 700 035, INDIA \\ (e-mail~: pradipta@isical.ernet.in)}
\date{}

\begin{document}

\maketitle
\vfill

\begin{abstract}
In this paper, we show that the projective tensor product of a
two-dimen\-sional $\ell^p$ space with a two-dimensional $\ell^q$
space never has the Mazur Intersection Property for a large
range of values of $p$ and $q$. For this purpose, we
characterise the extreme contractions from $\ell^p_2$ to
$\ell^q_2$ and obtain their closure.
\end{abstract}
\vfill

\hrule
\noindent {\bf Keywords and Phrases~:} Mazur Intersection Property, 
Extreme Contractions, $\ell^p$ spaces, Projective Tensor Product.\\
{\bf AMS Subject Classification~:} 46B20
\pagebreak

\section{Introduction}

A Banach space is said to have the Mazur Intersection Property
(MIP) if every closed bounded convex set is the intersection of
closed balls. In a finite-dimensional space $X$, this is
equivalent to the extreme points of the dual unit ball $B(X^*)$
being norm dense in the dual sphere $S(X^*)$. And, in general,

\begin{thm} \label{MT} For a Banach space $X$, the following are
equivalent~:
\bla
\item The w*-denting points of $B(X^*)$ are norm dense in
$S(X^*)$.
\item $X$ has the MIP.
\item Every support mapping on $X$ maps norm dense subsets of
$S(X)$ to norm dense subsets of $S(X^*)$.
\el 
\end{thm}
(see \cite{B2} or \cite{GGS} for details and related results)

Using this characterisation, Ruess and Stegall \cite{RS} have
shown that the injective tensor product of two Banach spaces of
dimension $\geq 2$ never has the MIP. And Sersouri \cite{S3} has
shown that in fact there is a two-dimensional compact convex set
in $X\otimes_{\e} Y$ that is not an intersection of balls.

The situation appears to be much more difficult for projective
tensor product spaces, since the extremal structure of the unit
ball of the dual of $X \otimes_{\p} Y$, i.e., ${\cal L} (X,
Y^*)$ (see e.g., \cite[Chapter~VIII]{DU}), is known only in some
very special cases and no pattern is discernible even in these
cases for a reasonable conjecture to be made in general. See
\cite{Grz5} or \cite{CHK} for a survey.

The simplest situation arises in a Hilbert space, where the
extreme contractions are characterised as isometries and
coisometries, by Kadison \cite{K} in the complex case (see also
\cite{H}) and by Grzaslewicz \cite{Grz2} in the real case. And
it immediately follows that the projective tensor product of two
Hilbert spaces never has the MIP.

Complications already increase significantly if we move on to
$\ell^p$-spaces. In fact, the complete picture eludes us even
for two-dimensional $\ell^p$-spaces. Here we show that the
projective tensor product of a two-dimensional $\ell^p$ space
with a two-dimensional $\ell^q$ space never has the MIP for a
large range of values of $p$ and $q$. For this purpose, we
characterise the extreme contractions from $\ell^p_2$ to
$\ell^q_2$ and obtain their closure. Some of the results about
extreme contractions were proved earlier in
\cite{Grz1,Grz2,Grz3} through different techniques. Our approach
is similar to that of \cite{CHK} for the case $p = q$ with
complex scalars. We, however, work only with real scalars. This
technique also lends itself naturally to the computation of the
closure.

A major portion of this work is contained in the first-named
author's Ph.\ D.\ Thesis \cite{PBTH} written under the
supervision of the second author. We take this opportunity to
thank a referee whose detailed comments on the paper led to
considerable improvement of its exposition.

\section{Extreme Contractions in ${\cal L}(\ell^p_2, \ell^q_2)$} 

{\sc Notations~:} For $1 < p < \iy$, $\bm{x} = (x_1, x_2) \in
\ell^p_2$ with $\|\bm{x}\| = 1$, define $\bm{x}^{p - 1} =
(sgn(x_1) |x_1|^{p - 1}$, $sgn(x_2) |x_2|^{p - 1})$ and
$\bm{x}^o = ( - x_2, x_1)$. Notice that, in general, $\bm{x}^{p
- 1}$ is the unique norming functional of \bm{x} and $\{\bm{x},
(\bm{x}^o)^{p - 1} \}$ is a basis for $\ell^p_2$, and if $p =
2$, $\bm{x}^{p - 1} = \bm{x}$ and $\{\bm{x}, \bm{x}^o
\}$ is orthonormal. Denote the vectors (1, 0) and (0, 1) by
$\bm{e}_1$ and $\bm{e}_2$ respectively.

We will need the following inequality \cite[Lemma 1.e.14]{LTz}
\begin{lem} \label{pq} Let $1 < p \leq q < \iy$. Let $\a = \pm
\sqrt{(p - 1)/(q - 1)}$. Then
\[ \left[ \rec{2} \left\{ \left| 1 + \a r \right|^q + \left| 1 -
\a r \right|^q \right\} \right]^{1/q} \leq \left[ \rec{2}
\left\{ \left| 1 + r \right|^p + \left| 1 - r \right|^p \right\}
\right]^{1/p} \] 
for all $r \in \R$ with strict inequality holding for $r \neq
0$. \end{lem}

\begin{thm} \label{tensor} For $1 < p$, $q < \iy$, an operator
$T : \ell^p_2 \lgra \ell^q_2$ with $\|T\| = 1$ is an extreme
contraction
\blr
\item {\rm \cite{Grz2}} for $p = q = 2$, if and only if $T$ is
an isometry.
\item {\rm \cite{Grz3}} for $p = 2 \neq q$, if and only if $T$
satisfies one of the following
\bla
\item $T$ attains its norm on two linearly independent vectors.
\item $T$ is of the form
\[ T = \alt{\bm{x} \otimes \bm{e}_i}{q < 2}{\bm{x} \otimes
\bm{y} + s \bm{x}^o \otimes (\bm{y}^o)^{q - 1}}{q > 2}\]
where $\bm{x}$ is any unit vector and, in the second case,
$|y_i|^q = \rec{2}$ and $s = \pm \rec{\sqrt{(q - 1)}} 2^{(q -
2)/q}$.
\el
\item {\rm \cite{Grz3}} for $p \neq 2 = q$, if and only if $T$
satisfies one of the following
\bla
\item $T$ attains its norm on two linearly independent vectors.
\item $T$ is of the form
\[ T = \alt{\bm{e}_i \otimes \bm{y}}{p > 2}{\bm{x}^{p - 1}
\otimes \bm{y} + s \bm{x}^o \otimes \bm{y}^o}{p < 2}\]
where $\bm{y}$ is any unit vector and, in the second case,
$|x_i|^p = \rec{2}$ and $s = \pm \sqrt{(p - 1)} 2^{(2 - p)/p}$.
\el
\item {\rm \cite{Grz1}} for $p = q \neq 2$, if and only if $T$
satisfies one of the following
\bla
\item $T$ attains its norm on two linearly independent vectors.
\item $T$ is of the form
\[T = \alt{\bm{e}_i \otimes \bm{y}}{p > 2, \qquad y_1 y_2 \neq
0}{\bm{x}^{p - 1} \otimes \bm{e}_j}{p < 2, \qquad x_1 x_2 \neq
0}\]
\el
\item for $1 < q < 2 < p < \iy$, if and only if $T$ satisfies
one of the following
\bla
\item $T$ attains its norm on two linearly independent vectors.
\item $T = \bm{x}^{p - 1} \otimes \bm{y}$ with \bm{x}, \bm{y}
unit vectors and $x_1 x_2 y_1 y_2 = 0$.
\el
\el \end{thm} 

\noindent {\sc Proof}. Let $T : \ell^p_2 \lgra \ell^q_2$,
$\|T\| = 1$. Then there exists $\bm{x} = (x_1, x_2) \in
\ell^p_2$ such that $\|\bm{x}\| = 1 = \|T\bm{x}\|$. Let $T\bm{x}
= \bm{y} = (y_1, y_2)$. Let $I_{xy} = \{ T : \|T\| \leq 1,
T\bm{x} = \bm{y} \}$. Then for any $T\in I_{xy}$, $(T -
\bm{x}^{p - 1} \otimes \bm{y})$ annihilates \bm{x} and so is of
rank $\leq 1$, whence $(T - \bm{x}^{p - 1} \otimes \bm{y}) =
\bm{x}^o \otimes \bm{u}$, for some $\bm{u}\in \ell^p_2$.
Further, $T^*(\bm{y}^{q - 1}) = \bm{x}^{p - 1}$, that is $(T^* -
\bm{y}\otimes \bm{x}^{p - 1})$ annihilates $\bm{y}^{q - 1}$,
whence $(T^* - \bm{y}\otimes \bm{x}^{p - 1}) = (\bm{y}^o)^{q -
1} \otimes \bm{v}$, for some $\bm{v}\in \ell^q_2$. Combining,
$T$ must be of the form
\[ T_s = \bm{x}^{p - 1} \otimes \bm{y} + s \bm{x}^o \otimes
(\bm{y}^o)^{q - 1}, \mbox{ for some } s \in \R \] 
In other words, $I_{xy} = \{ T_s : s \in \R, \|T_s\| \leq 1 \}$.
That is, $I_{xy}$ is a line segment (could be degenerate) in the
unit ball, and, its end points are extreme.

As in \cite{CHK}, pre- or post-multiplying by diag$(sgn(x_1),
sgn(x_2))$, diag$(sgn(y_1)$, $sgn(y_2))$ and permutation
matrices, if necessary --- each of which is an isometry --- we
may assume $x_1 \geq x_2 \geq 0$, $y_1 \geq y_2 \geq 0$.

For $r \in \R$, denote by $f_p (\bm{x}, r) = \bm{x} + r
(\bm{x}^o)^{p - 1}$ and $F_p (\bm{x}, r) = \|f_p (\bm{x},
r)\|^p$.

Then $F_p (\bm{x}, r) = |x_1 - r x_2^{p - 1}|^p + |x_2 + r
x_1^{p - 1}|^p$. Clearly, if $p = 2$ or $x_2 = 0$, $F_p (\bm{x},
r) = 1 + |r|^p$. Otherwise,
\beqa 
F_p (\bm{x}, r) 
& = & x_2^{ - p}\cdot |r x_2^p - x_1 x_2 |^p + x_1^{ - p}\cdot
|r x_1^p + x_1 x_2 |^p \\ 
& = & x_2^{ - p} G(r x_2^p - x_1 x_2) + x_1^{ - p} G(r x_1^p +
x_1 x_2)
\eeqa
where $G(u) = |u|^p$. Thus
\beqn 
\f{\partial}{\partial r} F_p (\bm{x}, r) & = & G'(r x_2^p - x_1
x_2) + G'(r x_1^p + x_1 x_2) \label{A1}
\eeqn
and $G'(0) = 0$, $G'(u) = p \cdot sgn(u) \cdot |u|^{p - 1}$.
Clearly, (\ref{A1}) also holds for $p = 2$ and $x_2 = 0$.

Now, $G'(u)$ is an odd function, positive and strictly
increasing for $u > 0$. Since the two arguments of $G'$ in
(\ref{A1}) add up to $r$, we have if $r > 0$ (resp.\ $r < 0$),
the one larger in absolute value is positive (resp.\ negative),
and so, $\f{\partial}{\partial r} F_p (\bm{x}, r)$ is positive
(resp.\ negative), i.e., $F_p (\bm{x}, r)$ is strictly
increasing (resp.\ decreasing) in $r > 0$ (resp.\ $r < 0$).

Further, if $p \neq 2$ and $x_2 \neq 0$
\beqn 
\left. \begin{array}{rcl} 
F_p (\bm{x}, 0) & = & 1, \qquad F'_p (\bm{x}, 0)~~ = ~~0, \\
F''_p (\bm{x}, 0) & = & p(p - 1)(x_1 x_2)^{p - 2}, \\ 
F'''_p (\bm{x}, 0) & = & p(p - 1)(p - 2)(x_1 x_2)^{p - 3} [x_1^p
- x_2^p], \\ 
F''''_p (\bm{x}, 0) & = & p(p - 1)(p - 2)(p - 3) (x_1 x_2)^{p -
4} [1 - 3(x_1 x_2)^p]
\end{array} \right\} \label{E2}
\eeqn
and so, if $1 < q < \iy$, for $H_{pq}(\bm{x}, r) = [F_p (\bm{x},
r)]^{q/p}$, we have,
\beqn 
\left. \begin{array}{rcl}
H_{pq}(\bm{x}, 0) & = & 1, \qquad H'_{pq}(\bm{x}, 0)~~ = ~~0, \\
H''_{pq}(\bm{x}, 0) & = & q(p - 1)(x_1 x_2)^{p - 2}, \\
H'''_{pq}(\bm{x}, 0) & = & q(p - 1)(p - 2)(x_1 x_2)^{p - 3}
[x_1^p - x_2^p], \\
H''''_{pq}(\bm{x}, 0) & = & q(p - 1)(x_1 x_2)^{p - 4} [(p - 2)(p
- 3) \\ 
&& \quad {} - 3(x_1 x_2)^p \{(p - 2)(p - 3) + (p - q)(p - 1)\}]
\end{array} \right\}
\label{E3}
\eeqn
where the derivatives are taken with respect to $r$.

Now, $T_s (f_p(\bm{x}, r)) = f_q (\bm{y}, rs)$, and thus for any
$r \neq 0$, $\|T_s (f_p(\bm{x}, r))\|^q = F_q (\bm{y}, rs)$ is
strictly increasing in $s \geq 0$ and strictly decreasing in $s
\leq 0$, and $F_q(\bm{y}, rs)$ is unbounded in $s$. Now, if $r
\neq 0$,
\[F_q(\bm{y}, 0) = 1 = [F_p(\bm{x}, 0)]^{q/p} < [F_p(\bm{x},
r)]^{q/p}\] 
So, there exists unique $s_ + (\bm{x}, \bm{y}, r) > 0$ and
unique $s_-(\bm{x}, \bm{y}, r) < 0$ such that
\beqn 
F_q(\bm{y}, rs_{\pm}) = [F_p(\bm{x}, r)]^{q/p} \label{E4}
\eeqn
And the quantity on the LHS becomes smaller or larger than the
one in the RHS according as $|s|$ gets smaller or larger.
Evidently, such $s_{\pm}$ also exist for $(\bm{x}^o)^{p - 1}$,
which we denote by $f_p(\bm{x}, \iy)$. In fact, in this case,
$|s_{\pm} (\bm{x}, \bm{y}, \iy)| = \|(\bm{x}^o)^{p -
1}\|/\|(\bm{y}^o)^{q - 1}\|$. Notice that $s_{\pm}$ is a
continuous function of $r \neq 0$ and elementary examples show
that $\lim_{r \ra 0} s_{\pm} (\bm{x}, \bm{y}, r)$ may not even
exist. Let
\[s^*_ + (\bm{x}, \bm{y}) = \inf \{s_ + (\bm{x}, \bm{y}, r) : r
\neq 0\} \quad s^*_ - (\bm{x}, \bm{y}) = \sup \{s_ - (\bm{x},
\bm{y}, r) : r \neq 0\}\]
Clearly, $T_s \in I_{xy}$ if and only if $s^*_ - \leq s \leq
s^*_ + $, i.e., $T_{s^*_{\pm}}$ are end points of $I_{xy}$ and
hence are extreme. Also let
\[\begin{array}{rclcl} 
s^{**}_ + (\bm{x}, \bm{y}) & = & \ds \liminf_{r\ra0} s_ +
(\bm{x}, \bm{y}, r) & \stackrel{\rm def}{=} & \ds \sup_{\e > 0}
\inf \{s_ + (\bm{x}, \bm{y}, r) : |r| < \e\} \\ 
s^{**}_ - (\bm{x}, \bm{y}) & = & \ds \limsup_{r\ra0} s_ -
(\bm{x}, \bm{y}, r) & \stackrel{\rm def}{=} & \ds \inf_{\e > 0}
\sup \{s_ - (\bm{x}, \bm{y}, r) : |r| < \e\}
\end{array}\]
Note that if we put $J_{xy} = \{ s : T_s$ is contractive in a
neighbourhood of \bm{x}\}, then $s^{**}_ - = \inf J_{xy}$ and
$s^{**}_ + = \sup J_{xy}$, though $s^{**}_{\pm}$ may not
necessarily belong to $J_{xy}$. Clearly, $s^{**}_ - \leq s^*_ -
\leq 0 \leq s^*_ + \leq s^{**}_ + $.

Now, either $s^*_{\pm}$ equals $s_{\pm}(\bm{x}, \bm{y}, r)$ for
some $r \neq 0$ (including $r = \iy$), in which case $s^*_{\pm}
\neq 0$ and $T_{s^*_{\pm}}$ attain their norm on two linearly
independent vectors, or $s^*_{\pm} = s^{**}_{\pm}$. 

Note that $T$ attains its norm on two linearly independent
vectors if and only if $T^*$ attains its norm on two linearly
independent vectors. Moreover, any such $T$ is exposed, and
hence strongly exposed.

Thus to complete the proof, the task that remains is to identify
all (if any) extreme contractions that attain their norm in
exactly one direction (called `of the desired type' in the
sequel). Then $s^*_{\pm} = s^{**}_{\pm}$, in which case $(i)$
$|s^{**}_{\pm}| < \iy$, $(ii)$ $s^{**}_{\pm} \in J_{xy}$, in
fact, $(iii)$ $T_{s^{**}_{\pm}} \in I_{xy}$.

Therefore, in different cases, we proceed to successively check
these three conditions and whenever we reach a contradiction, we
conclude that $s^*_{\pm} \neq s^{**}_{\pm}$ and $T_{s^*_{\pm}}$
is not of the desired type. And in case all the three conditions
are satisfied, we check whether it attains its norm in any
direction other than that of \bm{x} and only if it does not, we
get an extreme contraction of the desired type. This line of
reasoning is exemplified in the analysis of cases (II) and (IV)
below. However, in case (I), we can directly calculate
$s^*_{\pm}$.

{\sc Case(I) :} $(i)$ $p = 2$ and either $q = 2$ or $y_2 = 0 $;
$(ii)$ $q = 2$ and either $p = 2$ or $x_2 = 0 $; $(iii)$ $p \neq
2 \neq q \mbox{ and } x_2 = 0 = y_2$.
\beqa
T_s \mbox{ is a contraction } & \Lglra & F_q(\bm{y}, rs) \leq
[F_p(\bm{x}, r)]^{q/p} \quad \mbox{ for all } r \\ & \Lglra & 1
+ |rs|^q \leq [ 1 + |r|^p]^{q/p} \quad \mbox{ for all } r \\ &
\Lglra & |s|^q \leq \f{[ 1 + |r|^p]^{q/p} - 1}{|r|^q} \quad
\mbox{ for all } r \neq 0 
\eeqa

Note that the RHS $\equiv 1$ if $p = q2$ and is strictly
decreasing (resp.\ increasing) in $|r|$ for $q > p$ (resp.\ $q <
p$).

So, if $p = q$, $s_{\pm} \equiv \pm 1$, and hence, $s^*_{\pm} =
s^{**}_{\pm} = \pm 1$ and $T_{s^*_{\pm}}$ are isometries. And,
if $p \neq q$, the infimum of the RHS over $r
\neq 0$ yields
\[|s^*_{\pm}|^q = \alt{1}{q > p}{0}{q < p}\]

So, if $q < p$, $s^*_{\pm} = 0$, and $T_0$ is an extreme
contraction of the desired type. And if $q > p$, $s^*_{\pm} =
s_{\pm}(\iy) = \pm 1$ with $T_{\pm 1}$ attaining its norm at
both \bm{x} and $(\bm{x}^o)^{p - 1}$. It is interesting to note
that if $p \neq 2 \neq q$, $T_1$ in this case is the identity
operator.

For the remaining cases, we calculate $s^{**}_{\pm}$. Let
$\{r_n\}$ be a sequence of real numbers such that $r_n \lgra 0$
and $s_{\pm}(r_n) \lgra s^{**}_{\pm}$. If we assume
$|s^{**}_{\pm}| < \iy$, then $\{s_{\pm}(r_n)\}$ is a bounded
sequence. Now, by (\ref{E4}),
\beqn 
F_q(\bm{y}, r_n s_{\pm}(r_n)) = [F_p(\bm{x}, r_n)]^{q/p}
\label{E5}
\eeqn

{\sc Case (II) :} $q \neq 2$, $y_2 > 0$ and either $p = 2$ or
$x_2 = 0$.

In this case, subtracting 1 from both side of (\ref{E5}),
dividing by $r^2_n$ and taking limit as $n \lgra \iy$, we get by
L'Hospital's rule and (\ref{E2}) that $\mbox{LHS} \lgra
\rec{2} q (q - 1) s^2 (y_1 y_2)^{q - 2}$, where $s =
s^{**}_{\pm}$ and
\[\mbox{RHS} \lgra \alter{0}{p > 2}{\iy}{p < 2}{\ds \f{q}{2}}{p
= 2}\]

So, if $p < 2$, we have a contradiction, whence $T_{s^*_{\pm}}$
is not of the desired type, and if $p > 2$, $s^*_{\pm} =
s^{**}_{\pm} = 0$, and $T_0$ is extreme, and clearly of the
desired type.

If $p = 2$, we have
\beqn 
s^2(q - 1)(y_1 y_2)^{q - 2} = 1 \label{E6}
\eeqn

Now, if the $s$ given by (\ref{E6}) belongs to $J_{xy}$, we must
have
\beqn 
F_q (\bm{y}, rs) \leq [F_p (\bm{x}, r)]^{q/p} \quad \mbox{ for
all small } r \neq 0 \label{E7} 
\eeqn

Comparing the Taylor expansion of the two sides around $r = 0$
(for the LHS use (\ref{E2})), we see that the coefficients of
$1$, $r$ and $r^2$ on both sides are equal, whence the
inequality (\ref{E7}) for small $r$ implies the corresponding
inequality for the coefficient of $r^3$ on both sides, which,
for $r > 0$ and $r < 0$, leads to the equality
\beqn 
\rec{6}s^3 q(q - 1)(q - 2)(y_1 y_2)^{q - 3}(y_1^q - y_2^q) = 0
\label{E8}
\eeqn

Combining equations (\ref{E6}) and (\ref{E8}), we have $y_1^q =
y_2^q = 1/2$, $s^2 = \rec{(q - 1)} 4^{(q - 2)/q}$. But again the
equality in (\ref{E8}) pushes the inequality down to the
coefficients of $r^4$, i.e.,
\beqa
\rec{8} q(q - 2) & \geq & \rec{24} s^4 q(q - 1)(q - 2)(q - 3)
(y_1 y_2)^{q - 4} [1 - 3(y_1 y_2)^q] \\
\mbox{or } \quad 3 (q - 2) & \geq & \f{(q - 2)(q - 3)}{(q - 1)}
\eeqa

Now for $q < 2$, this leads to a contradiction, so that
$T_{s^*_{\pm}}$ is not of the desired type. On the other hand,
by Lemma~\ref{pq} for $p = 2$ and $q > 2$, we have that $T_s$
with the above parameters is a contraction that attains its norm
only in the direction of \bm{x} and hence, is of the desired
type.

{\sc Case (III) :} $p \neq 2$, $x_2 > 0$ and either $q = 2$ or
$y_2 = 0$.

This situation is dual to case (II) above.

{\sc Case(IV) :} $p \neq 2 \neq q$ and $x_2 > 0$, $y_2 > 0$.

In this case too, subtracting 1 from both side of (\ref{E5}),
dividing by $r^2_n$ and taking limit as $n \lgra \iy$, we get by
(\ref{E2}) and (\ref{E3})
\beqn 
(q - 1)(y_1y_2)^{q - 2} s^2 = (p - 1)(x_1x_2)^{p - 2} \label{E9}
\eeqn
where $s = s^{**}_{\pm}$. 

So, if $s \in J_{xy}$, comparing the Taylor expansion of the two
sides of (\ref{E7}) around $r = 0$ (use (\ref{E2}) for the LHS
and (\ref{E3}) for the RHS), by arguments similar to Case (II)
($p=2$), we must have
\beqn 
s^3(q - 1)(q - 2)(y_1y_2)^{q - 3}(y_1^q - y_2^q) = (p - 1)(p -
2)(x_1x_2)^{p - 3} (x_1^p - x_2^p) \label{E10}
\eeqn
and
\beqn 
s^4 (q - 1)(q - 2)(q - 3)(y_1 y_2)^{q - 4} [ 1 - 3(y_1 y_2)^q]
\leq (p - 1)(x_1 x_2)^{p - 4} \cdot \nonumber \\ 
{}[ (p - 2)(p - 3)\{ 1 - 3(x_1 x_2)^p \} - 3(p - q) (p - 1)(x_1
x_2)^p ]
\label{E11} 
\eeqn

Eliminating $s$ from (\ref{E9}) and (\ref{E10}) and using the
fact that \bm{x}, \bm{y} are unit vectors, we get
\beqn 
\f{(q - 2)^2}{(q - 1)}\left[\rec{(y_1y_2)^q} - 4 \right] = \f{(p
- 2)^2}{(p - 1)}\left[\rec{(x_1x_2)^p} - 4 \right] \label{E12}
\eeqn

Also, dividing (\ref{E11}) by the square of (\ref{E9}), we get
\beqn 
\f{(q - 2)(q - 3)}{(q - 1)} \left [\rec{(y_1 y_2)^q} - 3 \right
] \leq \f{(p - 2)(p - 3)}{(p - 1)} \left [ \rec{(x_1 x_2)^p} -
3\right] - 3 (p - q) \label{E13}
\eeqn

Notice that for $p = q$, we get from (\ref{E12}) that $x_i =
y_i$, $i = 1$, $2$, whence from (\ref{E9}), $s = \pm 1$, and
from (\ref{E10}), $y_1 = y_2 = x_1 = x_2$ for $s = - 1$. Thus,
\[ T_1 = \left[ \begin{array}{cc} 1 & 0 \\ 0 & 1 \end{array}
\right] \quad \mbox{ and } \quad T_{ - 1} = \left[
\begin{array}{cc} 0 & 1 \\ 1 & 0 \end{array} \right] \]
which, clearly, are isometries.

Now, let $p \neq q$. From (\ref{E12}) and (\ref{E13}), writing
$\ds \rec{(x_1 x_2)^p} = A$, we get 
\[\f{(p - 2)^2 (q - 3)}{(q - 2)(p - 1)}(A - 4) + \f{(q - 2)(q -
3)}{(q - 1)} \leq \f{(p - 2)(p - 3)}{(p - 1)}(A - 3) - 3 (p - q)
\]

Simplifying we get
\[\f{(q - p)(p - 2)}{(p - 1)(q - 2)} \cdot A \leq \f{2q(q - p)
(pq - p - q)}{(p - 1)(q - 1)(q - 2)}\]

So, if (a) $1 < q < 2 < p < \iy$, or, (b) $1 < p < q < 2$, or,
(c) $2 < p < q < \iy$, we have
\beqn 
A \leq \f{2q(pq - p - q)}{(p - 2)(q - 1)} \qquad \mbox{ i.e., }
\qquad A - 4 \leq \f{2(q - 2)(pq - p - q + 2)}{(p - 2)(q - 1)}
\label{E14}
\eeqn

And if (d) $1 < p < 2 < q < \iy$, or, (e) $1 < q < p < 2$, or,
(f) $2 < q < p < \iy$, we have
\beqn 
A \geq \f{2q(pq - p - q)}{(p - 2)(q - 1)} \qquad \mbox{ i.e., }
\qquad A - 4 \geq \f{2(q - 2)(pq - p - q + 2)}{(p - 2)(q - 1)}
\label{E15}
\eeqn

Now, $(x_1 x_2)^p = 1/A$ and $x_1^p + x_2^p = 1$, so $0 < 1/A
\leq 1/4$, i.e., $A \geq 4$ or $A - 4 \geq 0$.

But since $(pq - p - q + 2) = (p - 1)(q - 1) + 1$ is always
positive for $1 < p$, $q < \iy$, if $1 < q < 2 < p < \iy$, i.e.,
in case (a) above, we reach a contradiction at this point,
whence $T_{s^*_{\pm}}$ is not of the desired type. \qed

\noindent {\sc Remark}. $(a)$ For $p > q \geq 2$ and $y_2 = 0$,
the same result, as in Case (III) above, has been obtained by
Kan \cite[Lemma 6.2]{K2} for complex scalars too.

$(b)$ Recently we have come to know that P.\ Scherwentke
\cite{Sch} has proved a special case of Theorem~\ref{tensor}
$(v)$, i.e., when $p > 2$ and $1/p+1/q=1$, using techniques
similar to \cite{Grz3}.

\section{Partial Results in Remaining Cases}

In the last part of the proof of Theorem~\ref{tensor}, the
conditions (b) and (e) are dual to (c) and (f) respectively. And
in the cases (b) and (c), the inequality (\ref{E14}) implies
\beqn 
\rec{2} \leq x_1^p \leq \rec{2} \left[ 1 + \sqrt{\f{(q - 2)(pq -
p - q + 2)} {q(pq - p - q)}} \right] \label{E16}
\eeqn
while in cases (e) and (f), the inequality (\ref{E15}) implies
\beqn 
\rec{2} \left[ 1 + \sqrt{\f{(q - 2)(pq - p - q + 2)} {q(pq - p -
q)}} \right] \leq x_1^p < 1 \label{E17}
\eeqn

Now, in cases (b), (c), (e) and (f), we have from (\ref{E10})
that for $s < 0$, both sides of (\ref{E10}) must be 0, i.e., we
must have $x_1^p = x_2^p = 1/2 = y_1^q = y_2^q$. But then in
cases (e) and (f), we have a contradiction. So, in these two
cases, $s^*_ - $ gives extreme contractions not of the desired
type.

Also in case (d), (\ref{E15}) is always satisfied and
(\ref{E10}) implies that for $s > 0$, $x_1^p = x_2^p = 1/2 =
y_1^q = y_2^q$.

Now, from (\ref{E9}) it follows that $T_{s^{**}_{\pm}}$ is a
contraction if and only if
\beqn 
&& \left[ y_1^q \cdot \left| 1 + \a ( y_2/y_1 )^{q/2} r
\right|^q + y_2^q \cdot \left| 1 - \a ( y_1/y_2 )^{q/2} r
\right|^q \right]^{1/q} \nonumber\\ 
& \leq & \left[ x_1^p \cdot \left| 1 + ( x_2/x_1 )^{p/2} r
\right|^p + x_2^p \cdot \left| 1 - ( x_1/x_2 )^{p/2} r \right|^p
\right]^{1/p} \label{E18}
\eeqn
for all $r \in \R$, where $\a = \pm \sqrt{(p - 1)/(q - 1)}$ with
the sign being that of $s^{**}_{\pm}$. Thus for the particular
case of $x_1^p = x_2^p = 1/2 = y_1^q = y_2^q$, we have by
Lemma~\ref{pq} that in cases (b), (c) and (d) for both $s > 0$
and $s < 0$, we get extreme contractions of the desired type.

Thus, {\em modulo} duality, we are left with the following cases
unsolved :
\bln
\item Case (b) with $x_1^p > 1/2$ satisfying (\ref{E16}) for $s
> 0$ with $y_1$ given by (\ref{E12}).
\item Case (e) with $x_1$ satisfying (\ref{E17}) for $s > 0$
with $y_1$ given by (\ref{E12}).
\item Case (d) with $x_1^p > 1/2$ and $s < 0$ with $y_1$ given
by (\ref{E12}).
\el

In the remaining part of this section, we prove that in case
(b), i.e., for $1 < p < q < 2$, for $1/2 < x_1^p \leq 1/q$, we
get extreme contractions of the desired type. Specifically, we
prove 
\begin{lem} \label{pq1} Let $1 < p < q < 2$, $1/2 < x_1^p \leq
1/q$, then {\rm (\ref{E18})} holds for all $r \in \R$, with
equality only for $r = 0$.
\end{lem}

\noindent {\sc Proof}. For notational simplicity, put $x_i^p =
a_i$, $y_i^q = b_i$, $i = 1, 2$ and $a_2/a_1 = u$, $b_2/b_1 =
v$. Notice that, in this notation, (\ref{E12}) becomes
\beqn
(2 - q) \a v^{ - 1/2} (1 - v) & = & (2 - p) u^{ - 1/2} (1 - u)
\label{E19}
\eeqn
which implies $0 < v \leq u \leq 1$. It is also not difficult to
see that
\beqn
\a v^{ - 1/2} \leq u^{ - 1/2} & \Lglra & a_1 \leq \rec{q}
\label{E20}
\eeqn 
Also, in our notation (\ref{E18}) becomes
\beqn
&& \left[ \rec{1 + v} \cdot \left|1 + \a v^{1/2} r\right|^q +
\f{v}{1 + v} \cdot
\left|1 - \a v^{{} - 1/2} r\right|^q \right]^{1/q} \nonumber \\
& \leq & \left[ \rec{1 + u} \cdot \left|1 + u^{1/2} r\right|^p +
\f{u}{1 + u} \cdot
\left|1 - u^{{} - 1/2} r\right|^p \right]^{1/p} \label{D5}
\eeqn

{\sc Case I :} $0 \leq r \leq u^{1/2}$

Expanding $\mbox{LHS}^q$ and $\mbox{RHS}^p$ by Binomial series,
noting that $\a^2 = (p - 1)/(q - 1)$ and $(1 + x)^{q/p} \geq 1 +
\f{q}{p} x$ for all $x \geq 0$, it suffices to show that
\beqn
0 & \leq & (2 - q) \cdots (k - 1 - q) \a^{k - 2} v^{{} - (k -
2)/2} \left[\f{1 + ( - 1)^{k - 2} v^{k - 1}}{1 + v}\right]
\nonumber\\ & \leq & (2 - p) \cdots (k - 1 - p) u^{{} - (k -
2)/2} \left[\f{1 + ( - 1)^{k - 2} u^{k - 1}}{1 + u}\right]
\label{C6}
\eeqn
for all $k \geq 3$.

Now, as $0 < v < u < 1$, both $1 + ( - 1)^{k - 2} u^{k - 1}$ and
$1 + ( - 1)^{k - 2} v^{k - 1}$ are nonnegative, i.e., the first
inequality in (\ref{C6}) follows. Also for $k = 3$, the second
inequality in (\ref{C6}) is an equality by (\ref{E19}). And for
$k \geq 4$, dividing both sides of (\ref{C6}) by that of
(\ref{E19}), it suffices to show
\beqa
&& (3 - q) \cdots (k - 1 - q) \a^{k - 3} v^{{} - (k - 3)/2}
\left[\f{1 + ( - 1)^{k - 2} v^{k - 1}}{1 - v^2}\right] \\ & \leq
& (3 - p) \cdots (k - 1 - p) u^{{} - (k - 3)/2} \left[\f{1 + ( -
1)^{k - 2} u^{k - 1}}{1 - u^2}\right]
\eeqa

But for $k \geq 4$, $(3 - q) \cdots (k - 1 - q) < (3 - p) \cdots
(k - 1 - p)$, and by (\ref{E20}), $\a^{k - 3} v^{{} - (k - 3)/2}
\leq u^{{} - (k - 3)/2}$. Also, it is not difficult to see that
for any $k \geq 4$, $[1 + ( - 1)^{k - 2} x^{k - 1}]/(1 - x^2)$
is strictly increasing for $0 < x < 1$. Since $0 < v < u < 1$,
Case I follows.

Notice that for $r = u^{1/2}$, we get
\begin{equation}
\left[ \rec{1 + v} \cdot \left(1 + \a v^{1/2} u^{1/2}\right)^q +
\f{v}{1 + v} \cdot
\left(1 - \a v^{{} - 1/2} u^{1/2}\right)^q \right] \leq \left[
(1 + u)^{p - 1}
\right]^{q/p} \label{B3}
\end{equation}

{\sc Case II :} $r \leq 0$ and $r \geq u^{1/2}$.

Notice that if we put $t = - u^{ - 1/2} r/(1 - u^{ - 1/2} r)$,
(\ref{D5}) becomes
\beqn
&& \left[ \rec{1 + v} \cdot \left|1 - (1 + c) t\right|^q +
\f{v}{1 + v} \cdot \left|1 - (1 - c/v) t\right|^q \right]^{1/q}
\nonumber\\ & \leq & \left[ \rec{1 + u} \cdot \left|1 - (1 + u)
t\right|^p + \f{u}{1 + u}
\right]^{1/p} \label{B1} 
\eeqn
where $c = \a u^{1/2} v^{1/2}$, and the ranges $r \leq 0$ and $r
\geq u^{1/2}$ become $0 \leq t < 1$ and $t > 1$ respectively.
Thus, we have to prove (\ref{B1}) for $t \geq 0$, $t \neq 1$.

Notice that by (\ref{E20}), $c \leq v$ and $c \leq \a u < u$.
Put
\beqa
\ph_1(t) & = & \rec{1 + v} \left[ \left|1 - (1 + c) t\right|^q +
v \cdot |1 - (1 - c/v) t|^q \right] \mbox{ and}\\
\ph_2(t) & = & \rec{1 + u} \left[ \left|1 - (1 + u) t\right|^p +
u \right]
\eeqa

Put $f(t) = q \log \ph_2(t) - p \log \ph_1(t)$. We have to show
$f(t) > 0$ for $t \neq 0$. Now, $f'(t) = (q \ph_1(t) \ph_2'(t) -
p \ph_1'(t)
\ph_2(t))/(\ph_1(t) \ph_2(t))$, so that $f'(t) > , = , \mbox{ or
} < 0$ according as $q \ph_1(t) \ph_2'(t) - p \ph_1'(t) \ph_2(t)
> , = , \mbox{ or } < 0$; or, equivalently,
\beqn 
&& sgn[1 - (1 - c/v)t]|1 - (1 - c/v)t|^{q - 1} \cdot [1 - (uv +
c) \cdot g(t)]
\nonumber \\
& < , = , \mbox{ or } > & sgn[1 - (1 + c)t] \cdot |1 - (1 +
c)t|^{q - 1} \cdot [1 + (u - c) \cdot g(t)] \label{B2}
\eeqn
where $g(t) = \{1 - sgn[1 - (1 + u)t] \cdot |1 - (1 + u)t|^{p -
1}\}/c(1 + u)$.

Notice that $g'(t) = c^{ - 1} (p - 1) |1 - (1 + u) t|^{p - 2}
\geq 0$, and hence, $g(t)$ is strictly increasing with $g(0) =
0$.

{\sc Subcase 1 :} $0 \leq t \leq 1/(1 + c)$.

Since $f(0) = 0$, and it suffices to prove $f'(t) \geq 0$, or,
in (\ref{B2}), LHS $\leq $ RHS. Notice that in this case, every
factor on the two sides of (\ref{B2}), except possibly the third
term on the LHS, is nonnegative. And the third term on the LHS
is decreasing, positive at $t = 0$ and is $\leq 0$ at $t = 1/(1
+ u)$. When this term is $\leq 0$, we have nothing to prove. And
thus it suffices to prove
\[ \left[\f{1 - (1 - c/v)t}{1 - (1 + c)t} \right]^{q - 1} \leq
\f{[1 + (u - c)
\cdot g(t)]}{[1 - (uv + c) \cdot g(t)]} \]
for the values of $t$ for which $g(t) < 1/(uv + c)$, which
exclude the values $1/(1 + u) \leq t \leq 1/(1 + c)$.

Again since in this range all the factors are positive and the
two sides are equal at $t = 0$, taking logarithm and
differentiating, it suffices to show
\[ \f{(q - 1) c}{v [1 - (1 - c/v)t] \cdot [1 - (1 + c)t]} \leq
\f{u g'(t)}{[1 + (u - c) \cdot g(t)] \cdot [1 - (uv + c) \cdot
g(t)]}\]

Simplifying the expressions and putting $s = (1 + u)t$, this is
equivalent to
\[A \cdot [\f{c^2}{(1 - s)} + v (1 - s)^{p - 1}] + D \cdot [c^2
(1 - s) + u^2 v (1 - s)^{1 - p}] + B \cdot [c^2 - uv] \geq 0\]
where
\beqa 
A & = & (uv + c)(u - c) > 0\\ B & = & (uv + c)(1 + c) + (v -
c)(u - c)\\ D & = & (1 + c)(v - c) \geq 0
\eeqa

Now, in the range $0 \leq s < 1$ and so, we can expand the LHS
by Binomial and geometric series. Note that
\[ A + D + B = v (1 + u)^2 \mbox{ and } A + u^2 D - u B = - c^2
(1 + u)^2 \] whence the constant term on the LHS is
\[c^2(A + D + B) + v(A + u^2 D - u B) = 0\]

On the other hand, we have from (\ref{E19}), that
\beqa
A - D & = & (1 + u) [ c(1 - v) - v(1 - u) ] = \f{c(1 + u)(1 -
v)(q - p)}{2 - p}\\ A - u^2 D & = & c(1 + u) [ u(1 - v) - c(1 -
u)] = uv(1 - u^2) \f{q - p}{(q - 1)(2 - q)} > 0
\eeqa
whence the coefficient of $s$ on the LHS is
\[c^2(A - D) - (p - 1)v (A - u^2 D) = c^2 (1 + u)(q - p) \left[
\f{c(1 - v)}{2 - p} - \f{v(1 - u)}{2 - q} \right] = 0 \]

Therefore, and since $1 < p < 2$, we have the coefficient of
$s^k$, $k \geq 2$, is
\beqa 
&& A \left\{c^2 - \f{v(p - 1)(2 - p) \cdots (k - p)}{k!} \right\} +
D u^2 v \f{(p - 1)p(p + 1) \cdots (p + k - 2)}{k!} \\
&& \qquad \geq A c^2 - \f{(p - 1)(2 - p) \cdots (k - p)}{k!} \cdot
v (A - u^2 D) \\
&& \qquad = D c^2 + v (A - u^2 D) (p - 1) \left[1 - \f{(2 - p)
\cdots (k - p)}{k!} \right] > 0
\eeqa 

{\sc Subcase 2 :} $1/(1 + c) < t < v/(v - c)$.

Since $(|x| + |y|)^a \geq |x|^a + |y|^a$ for $a > 1$, we have
that in (\ref{B1})
\beqa
\mbox{RHS}^q & \geq & \left(\rec{1 + u}\right)^{q/p} [(1 + u)t -
1]^q + \left( \f{u}{1 + u} \right)^{q/p} \\
\mbox{LHS}^q & = & \rec{1 + v} \cdot [(1 + c) t - 1]^q + \f{v}{1
+ v} \cdot [1 - (1 - c/v) t]^q
\eeqa

Comparing the first term of the two sides, it suffices to show
\beqa 
\rec{1 + v} \cdot [(1 + c) t - 1]^q & \leq & \left(\rec{1 +
u}\right)^{q/p} [(1 + u)t - 1]^q \\
\mbox{or, } \left[ \f{(1 + c) t - 1}{(1 + u) t - 1} \right] &
\leq & \f{(1 + v)^{1/q}}{(1 + u)^{1/p}}
\eeqa
for $t > 1/(1 + c)$, the LHS is increasing, and the maximum
value at ``$t =
\iy$" is $(1 + c)/(1 + u)$. Thus it suffices to show that 
\[ (1 + c) \leq (1 + v)^{1/q} \cdot (1 + u)^{1 - 1/p} \]
but this follows from (\ref{B3}).

And comparing the second term, we need to show
\[\left( \f{u}{1 + u} \right)^{q/p} \geq \f{v}{1 + v} \cdot [1 -
(1 - c/v) t]^q \]

For $1/(1 + c) < t < v/(v - c)$, the RHS is decreasing and it
suffices to prove
\[\left( \f{u}{1 + u} \right)^{q/p} \geq \f{v}{1 + v} \cdot
\left[1 - \f{(1 - c/v)}{(1 + c)} \right]^q = \f{v}{1 + v} \cdot
\left[\f{(c + c/v)}{(1 + c)} \right]^q
\]

Now, if we consider the function
\[ h(x) = \rec{1 + v} \cdot \left[\f{x - c}{1 + x}\right]^q +
\f{v}{1 + v} \cdot
\left[ \f{x + c/v}{1 + x}\right]^q \]
we see that
\[h'(x) = q \cdot \f{(1 + c) (x - c)^{q - 1} + (v - c) (x +
c/v)^{q - 1}}{(1 + v) (1 + x)^{q + 1}}\] whence $h(x)$ is
increasing for $x \geq c$, and the inequality (\ref{B1}) at $t =
1/(1 + u)$, yields
\[ \left( \f{u}{1 + u} \right)^{q/p} \geq h(u) \geq h(c) \]
This proves the subcase 2.

{\sc Subcase 3 :} $t \geq v/(v - c)$. Notice that if $v = c$,
this case does not arise.

In this case, the LHS of (\ref{B2}) is $\geq 0$, while the RHS $
< 0$, whence $f'(t) < 0$, and the minimum value of $f$ is
attained for ``$t = \iy$". Now that this value is $ > 0$ follows
from (\ref{B3}). This completes the proof of Case II, and hence,
of the Lemma. \qed

In the particular case $a_1 = 1/q$, replacing $r/\sqrt{q - 1}$
by $t$, we get the following interesting inequality, the case $q
= 2$ being immediate from Lemma~\ref{pq}.
\begin{cor} Let $1 < p < q \leq 2$. Then 
\[\left[\rec{p} |1 + (p - 1) t|^q + \f{p - 1}{p} |1 - t|^q
\right]^{1/q} \leq
\left[\rec{q} |1 + (q - 1) t|^p + \f{q - 1}{q} |1 - t|^p
\right]^{1/p} \] for all $t \in \R$ with strict inequality
holding for $t \neq 0$. \end{cor}

\noindent {\sc Remark}. We do not know whether the range $1/2 \leq x_1^p
\leq 1/q$ exhausts all values of $x_1^p$ for which we get a
contraction. However, it is not very difficult to see that we
cannot have the entire range in (\ref{E16}). Indeed, when $x_1^p$
is the right end point, we do not even get a contractive $T$, as
in that case tracing our arguments back we find that the
coefficients of $r^4$ in the Taylor expansion of the two sides
of (\ref{E7}) must be equal, and hence, as in the case of $r^3$,
we must have equality of the coefficients of $r^5$ as well. But
then direct computations reveal a contradiction.

\section{The Closure of Extreme Contractions}

We now obtain the closure of the extreme contractions in the
cases described in Theorem~\ref{tensor}.
\begin{thm} \label{close} In all the cases described in
Theorem~\ref{tensor}, except the case $p = q \neq 2$, the
set of extreme contractions is closed.

And in the case $p = q \neq 2$, the closure of extreme
contractions may have operators of the form diag$(1, s)$, $|s| <
1$ (upto isometric factors of signum or permutation matrices) as
additional elements.
\end{thm}

\noindent {\sc Proof}. In case $(i)$, the result is obvious. And
case $(iii)$ is dual to case $(ii)$.  Now, in the cases $(ii)$
and $(v)$, the set of operators of the form (b) is clearly
closed. And in case $(iv)$, the closure of the set of operators
of the form (b) contains only the operators $\bm{e}_i \otimes
\bm{e}_j$, $i$, $j = 1$, $2$ in addition.

Let us consider the set of operators of the type (a) in cases
$(ii)$, $(iv)$ and $(v)$. Let $\{T_n\}$ be a sequence of
operators of the type (a). Let $T_n \lgra T$ in operator norm.
Let $\bm{x}_n = (x_{n1}, x_{n2})$ be such that $\|\bm{x}_n\| = 1
= \|T_n\bm{x}_n\|$. Let $T_n\bm{x}_n = \bm{y}_n = (y_{n1},
y_{n2})$. Then $T_n$ is of the form
\[ T_n = \bm{x}_n^{p - 1} \otimes \bm{y}_n + s^*_{\pm}(\bm{x}_n,
\bm{y}_n) \bm{x}_n^o \otimes (\bm{y}_n^o)^{q - 1}\]
where $s^*_{\pm}(\bm{x}_n, \bm{y}_n)$ is as in our earlier
discussion. For notational simplicity, write $s^*_{\pm}
(\bm{x}_n, \bm{y}_n) = \pm s_n$. Passing to a subsequence, if
necessary, assume $\bm{x}_n \lgra \bm{x} = (x_1, x_2)$,
$\bm{y}_n \lgra \bm{y} = (y_1, y_2)$ (by compactness of the unit
balls of $\ell^p_2$ and $\ell^q_2$), and all the $s_n$'s have
the same sign, without loss of generality, positive.

Clearly, $\|T\| = 1$ and $T \bm{x} = \bm{y}$, whence $T$ is of
the form
\[ T = \bm{x}^{p - 1} \otimes \bm{y} + s \bm{x}^o \otimes
(\bm{y}^o)^{q - 1}\]

Also, as $T_n \lgra T$, $s_n = \|T_n - \bm{x}_n^{p - 1} \otimes
\bm{y}_n\|/ \|\bm{x}_n^o\| \cdot \|(\bm{y}_n^o)^{q - 1}\| \lgra
\|T - \bm{x}^{p - 1} \otimes \bm{y}\|/ \|\bm{x}^o\| \cdot
\|(\bm{y}^o)^{q - 1}\|$, i.e., $\{s_n\}$ is convergent. Clearly,
$s_n \lgra s$.

Now, since $T_n$ is of the type (a), there exists $\bm{z}_n =
\bm{x}_n + r_n (\bm{x}_n^o)^{p - 1}$ with \mbox{$r_n \neq 0 \in
\R$} such that $\|T_n \bm{z}_n\| = \|\bm{z}_n\|$. Again we may
assume all $r_n$'s are of the same sign, in particular positive
and $r_n \lgra r$, where $0 \leq r \leq \iy$. If $0 < r < \iy$,
$\bm{z}_n \lgra \bm{z} = \bm{x} + r (\bm{x}^o)^{p - 1}$ and
$\|T\bm{z}\| = \|\bm{z}\|$, i.e., $T$ is also of the type (a).
Also, if $r_n \lgra \iy$, let $\bm{u}_n = \bm{z}_n /
\|\bm{z}_n\|$. Then $\bm{u}_n \lgra \bm{u} = (\bm{x}^o)^{p - 1}
/ \|(\bm{x}^o)^{p - 1}\|$ and $\|T\bm{u}\| = 1$, so that $T$
again is of the type (a).

Now, suppose $r_n \lgra 0$. Then from $\|T_n\bm{z}_n\| =
\|\bm{z}_n\|$ we have
\beqn 
F_q(\bm{y}_n, r_n s_n) - [F_p(\bm{x}_n, r_n)]^{q/p} = 0
\label{A12}
\eeqn

For $(ii)$, if $p = 2$ and $q < 2$, since $T_n$ is of type (a),
we have $y_{n1} y_{n2} \neq 0$ for all $n$. And if $q > 2$, we
have two possibilities; either there is a subsequence for which
$y_{n1} y_{n2} = 0$, or, eventually $y_{n1} y_{n2} \neq 0$. In
the first case, we restrict ourselves only to that subsequence,
and we have, by case (I), $s_n = 1$ for all $n$, whence $s = 1$.
Also, $y_1 y_2 = 0$. So, $T = \bm{x} \otimes \bm{e}_i + \bm{x}^o
\otimes \bm{e}_j$ $(i \neq j)$, and it is clear that $T$ is of
type (a) (see case (I)). And in the second case, we assume
$y_{n1} y_{n2} \neq 0$ for all $n$. Then dividing (\ref{A12}) by
$r_n^2$ and taking limit as $n \lgra \iy$, we get by
L'Hospital's rule
\beqn 
(q - 1) s^2 (y_1 y_2)^{q - 2} - 1 = 0 \label{A13}
\eeqn

If $q > 2$, for $y_1 y_2 = 0$, this leads to a contradiction,
whence $y_1 y_2 \neq 0$. Then (\ref{E6}) and (\ref{A13})
coincides, i.e., we have $s = s^{**}_ + (\bm{x}, \bm{y})$. Now,
our analysis in case (II) shows that only for $y_1^q = y_2^q =
1/2$, $s^{**}_ + $ gives a contraction (which is an extreme
contraction of type (b)). And in every other case, we run into a
contradiction, i.e., we must have $r_n \not \ra 0$.

And if $q < 2$, for $y_1 y_2 = 0$, (\ref{A13}) makes sense only
if $s = 0$. In that case, $T = \bm{x} \otimes \bm{e}_i$, which,
by case (I), is an extreme contraction of the type (b). And for
$y_1 y_2 \neq 0$, we again have $s = s^{**}_ + (\bm{x}, \bm{y})$
and our analysis in case (II) shows that this case always leads
to a contradiction.

So, in both the cases, the closure of the set of operators of
the type (a) contains at most operators of type (b), and
therefore, the set of extreme contractions is closed.

In $(iv)$, i.e., if $p = q$, by duality, it suffices to consider
$p > 2$. Since $T_n$ is of type (a), we have three
possibilities; (1) either there is a subsequence for which both
$x_{n1} x_{n2} = 0$ and $y_{n1} y_{n2} = 0$, or, (2) there is a
subsequence for which $x_{n1} x_{n2} \neq 0$ and $y_{n1} y_{n2}
= 0$, or, (3) eventually both $x_{n1} x_{n2} \neq 0$ and $y_{n1}
y_{n2} \neq 0$.

In the first case, we again restrict ourselves only to that
subsequence, and we have, by case (I), $s_n = 1$ for all $n$,
whence $s = 1$. Also, $x_1 x_2 = 0$ and $y_1 y_2 = 0$. Now again
by case (I), $T$ is of type (a).

In cases (2) and (3), dividing (\ref{A12}) by $r_n^2$ and taking
limit --- through a subsequence if necessary --- as $n \lgra
\iy$, we get in the second case
\[(x_1 x_2)^{p - 2} = 0\] 
and in the third case
\[ (y_1y_2)^{p - 2} s^2 = (x_1x_2)^{p - 2} \]

So, in case (2), $y_1 y_2 = 0$, and we have a contradiction
unless $x_1 x_2 = 0$. And in that case, $T$ is of the form
diag$(1, s)$ upto isometric factors of signum or permutation
matrices. Now, $T$ is a contraction for $ - 1 \leq s \leq 1$ and
is extreme (in fact, an isometry) only for $s = \pm 1$. However,
we do not know precisely if they actually belong to the closure.

In case (3), if $x_1 x_2 = 0$, we get a contradiction unless
$y_1 y_2 = 0$ or $s = 0$. If $y_1 y_2 \neq 0$, $s = 0 =
s^*_{\pm} (\bm{x}, \bm{y})$, whence $T$ is extreme. And if $y_1
y_2 = 0$, we get the conclusions as in case (2). If $x_1 x_2
\neq 0$, $y_1y_2 = 0$ leads to a contradiction, and if $y_1 y_2
\neq 0$, $s = s^{**}_ + (\bm{x}, \bm{y})$, so that $T$ is an
isometry and hence is of type (a).

In case $(v)$, i.e., if $1 < q < 2 < p < \iy$, since $T_n$ is of
type (a), we must have $x_{n1} x_{n2} y_{n1} y_{n2} \neq 0$ for
all $n$. And a similar argument leads to
\[ s^2 (q - 1)(y_1y_2)^{q - 2} = (p - 1)(x_1x_2)^{p - 2} \]

If $x_1 x_2 \neq 0$ the only situation that does not lead to any
contradiction --- either immediate or to the fact that $T$ is a
contraction --- is both $y_1 y_2 = 0$ and $s = 0$. And in that
case, $s = 0 = s^*_{\pm}(\bm{x}, \bm{y})$, so that $T$ is
extreme. And if $x_1 x_2 = 0$, we must have $s = 0$, in which
case, by cases (I) and (II), $s = 0 = s^*_{\pm} (\bm{x},
\bm{y})$ and $T$ is extreme. Thus in this case too, the set of
extreme contractions is closed. \qed

\begin{thm} In each of the following cases of $1 < p$, $q <
\iy$, $\ell^p_2 \otimes_{\p} \ell^q_2$ lacks the MIP~:
\blr
\item $p$ and $q$ are conjugate exponents, i.e., $\rec{p} +
\rec{q} = 1$.
\item Either $p$ or $q$ is equal to $2$.
\item $2 < p$, $q < \iy$.
\el \end{thm}

\noindent {\sc Proof}. The dual of $\ell^p_2 \otimes_{\p}
\ell^q_2$ is ${\cal L}(\ell^p_2, \ell^{q'}_2)$, where $\rec{q} +
\rec{q'} = 1$ and the closure of extreme contractions in none of
the above cases contains norm 1 operators of the form $\bm{x}^{p
- 1} \otimes \bm{y}$, where \bm{x} and \bm{y} are unit vectors
with $x_1 x_2 y_1 y_2 \neq 0$. \qed

\noindent {\sc Remark}. The fact that operators of the above
form do not belong to the closure of extreme contractions in any
of these cases seems to suggest that this is a general
phenomenon. It is possible that this is happens in higher
dimensions as well. Can one give a proof of this without
precisely characterising the extreme contractions? What seems to
be required is a more tractable necessary condition for
extremality, or, for belonging to the closure of extreme
contractions.


\begin{thebibliography}{99}

\bibitem{PBTH} Bandyopadhyay P., {\em The Mazur Intersection
Property in Banach Spaces and Related Topics}, Ph.~D. Thesis,
Indian Statistical Institute, Calcutta (February, 1991).

\bibitem{B2} P.\ Bandyopadhyay, {\em The Mazur Intersection
Property for Families of Closed Bounded Convex Sets in Banach
Spaces}, Colloq.\ Math.\ {\bf LXIII} (1992), 45--56.

\bibitem{DU} Diestel J.\ and Uhl J.\ J., Jr., {\em Vector
Measures}, Math.\ Surveys {\bf 15}, Amer.\ Math.\ Soc.,
Providence, R.\ I.\ (1977).

\bibitem{GGS} J.\ R.\ Giles, D.\ A.\ Gregory and B.\ Sims, {\em
Characterisation of Normed Linear Spaces with the Mazur's
Intersection Property}, Bull.\ Austral.\ Math.\ Soc.\ {\bf 18}
(1978), 105--123.

\bibitem{Grz1} Grzaslewicz R., {\em Extreme Operators on
$2$-dimensional $\ell^p$-spaces}, Colloq.\ Math.\ {\bf 44}
(1981), 309--315.

\bibitem{Grz2} Grzaslewicz R., {\em Extreme Contractions on Real
Hilbert Spaces}, Math.\ Ann.\ {\bf 261} (1982), 463--466.

\bibitem{Grz4} Grzaslewicz R., {\em Exposed Points of the Unit
Ball of ${\cal L}(H)$}, Math.\ Z.\ {\bf 193} (1986), 595--596.

\bibitem{Grz3} Grzaslewicz R., {\em Extremal Structure of ${\cal
L}(\ell^2_m, \ell^p_n)$}, Linear and Multilinear Algebra {\bf
24} (1989), 117--125.

\bibitem{Grz5} Grzaslewicz R., {\em Survey of Main Results about
Extreme Operators on Classical Banach Spaces}, (preprint).

\bibitem{H} Halmos P.\ R., {\em A Hilbert Space Problem Book},
New Jersey (1967).

\bibitem{K} Kadison R.\ V., {\em Isometries of Operator
Algebras}, Ann.\ Math.\ {\bf 54} (1951), \mbox{325--338}.

\bibitem{CHK} Kan C.\ H., {\em A Class of Extreme $L_p$
Contractions, $p \neq 1$, $2$, $\iy$, and Real $2 \times 2$
Extreme Matrices}, Illinois J.\ Math., {\bf 30} (1986),
612--635.

\bibitem{K2} Kan C.\ H., {\em Norming Vectors of Linear
Operators between $L_p$ Spaces}, Pacific J.\ Math., {\bf 150}
(1991), 309--327.

\bibitem{LTz} Lindenstrauss J.\ and Tzafriri L., {\em Classical
Banach Spaces}, Vol.\ II, Ergebnisse der Mathematik und ihrer
Grenzgebiete, Band {\bf 97}, Springer Verlag (1979).

\bibitem{RS} Ruess W.\ M.\ and Stegall C.\ P., {\em
Weak*-Denting Points in Duals of Operator Spaces}, Proceedings
of the Missouri Conference on Banach Spaces, Springer Lecture
Notes {\bf 1166} (1985), 158--168.

\bibitem{Sch} Scherwentke P., {\em On the Geometry of ${\cal
L}(\ell^p_2, \ell^q_2)$ and $\ell^q_2 \otimes_{\e} \ell^q_2$},
Arch.\ Math., (to appear).

\bibitem{S3} Sersouri A., {\em Smoothness in Spaces of Compact
Operators}, Bull.\ Austral.\ Math.\ Soc.\ {\bf 38} (1988),
221--225.

\end{thebibliography}
\end{document}